\documentclass[11 pt,leqno]{amsart}

\usepackage[utf8]{inputenc} 
\usepackage{amsmath,amsthm,amssymb,amscd,graphicx}
\usepackage{color}
\usepackage[normalem]{ ulem }
\usepackage{soul}
\usepackage[left=2.5cm,right=2.5cm,top=3cm,bottom=3cm]{geometry}
\usepackage{hyperref}
\usepackage{enumerate}
\usepackage{mathrsfs}
\usepackage{multirow}

\newtheorem{thm}{Theorem}[section]
 
 \newtheorem{defi}[thm]{Definition}

 \newtheorem{prop}[thm]{Proposition}
 
 \theoremstyle{definition}
 
 \theoremstyle{remark}

 \numberwithin{equation}{section}
\def\be#1 {\begin{equation} \label{#1}}
\newcommand{\ee}{\end{equation}}
\renewcommand{\phi}{\varphi}

\def\s{\sigma}
\def\C{\mathbb C}
\def\R{\mathbb R}

\def\N{\mathbb N}

\def\H{\mathcal H}
\def\W{\mathcal W}

\def\e{e}

\def\eps{\epsilon}
\def\dis{\displaystyle}    

 \newcommand{\om}{  \omega   }

\newcommand{\ov}{  \overline  }
\newcommand{\p}{  {\bf p}  }

\renewcommand\>{\rangle}

\definecolor{gr}{rgb}   {0.,   0.69,   0.23 }
\definecolor{bl}{rgb}   {0.,   0.5,   1. }
\definecolor{mg}{rgb}   {0.85,  0.,    0.85}
\definecolor{yl}{rgb}   {0.8,  0.7,   0.}
\definecolor{or}{rgb}  {0.7,0.2,0.2}

\begin{document}

\author{ Laurent Thomann }
\address{Institut  \'Elie Cartan, Universit\'e de Lorraine, B.P. 70239,
F-54506 Vand\oe uvre-l\`es-Nancy Cedex}
\email{laurent.thomann@univ-lorraine.fr}

\author{Nicolas Burq}
\address{Laboratoire de Math\'ematiques d'Orsay,
CNRS, Universit\'e Paris--Saclay, B\^at. 307, F-91405 Orsay Cedex, and Institut Universitaire de France}
\email{nicolas.burq@universite-paris-saclay.fr}
 
 \thanks{N. Burq is  partially supported by the ANR projet Smooth "ANR-22-CE40-0017" and Institut Universitaire de France}
 
\thanks{L. Thomann  is  partially supported by the ANR projet Smooth "ANR-22-CE40-0017" and a Tohoku University-Universit\'e de Lorraine joint research fund}

\title[Evolution of Gaussian measures and  NLS]{Evolution of Gaussian measures and application to  the one dimensional nonlinear Schr\"odinger equation}

\begin{abstract}
In this note, we give an overview of some results obtained in \cite{Burq-Thomann}. This latter work is devoted to the study of the one-dimensional nonlinear Schr\"odinger equation with random initial conditions. Namely, we describe the nonlinear evolution of Gaussian measures and we deduce global well-posedness and scattering results for the corresponding nonlinear Schr\"odinger equation.
\end{abstract}

\maketitle
 

   \section{Introduction}

The motivation of this work is to study the long time behaviour of the nonlinear Schr\"odinger equation with random initial conditions
 \begin{equation*} 
  \left\{
      \begin{aligned}
         &i\partial_sU+\partial_{y}^2U=|U|^{p-1}U,\quad (s,y) \in  \R\times \R,
       \\  &  U(0)  =U_0 \in L^2(\R).
      \end{aligned}
    \right.
\end{equation*} 
 We will be able to prove: \medskip

 $\blacktriangleright$ almost sure global existence results ($p>1$)\medskip
 
  $\blacktriangleright$ almost sure scattering results ($p>3$). \medskip
 
 To prove these results, we will first construct measures on the space of initial data for which we can describe precisely the  non trivial evolution  by the linear Schr\"odinger flow. Then we prove that the  nonlinear evolution of these measures is absolutely continuous with respect to their  linear evolutions, with  quantitative estimates on the Radon-Nikodym derivative. To the best of our knowledge, these results  are the first ones giving insight,  in a non compact setting on the time evolution of the statistical distribution of solutions of a nonlinear PDE. They also are the first ones providing scattering for NLS for large initial data without assuming decay at infinity: our solutions are in a slighly larger space than   $L^2(\R)$.


 \subsection{On invariant measures}
 To begin with, let us recall the definition of a measure left  invariant   by a one parameter group.
\begin{defi}
Consider a  space $X$ and a one parameter group $(\Phi(t,.))_{t\in \R}$ with  
$$\Phi(t,.):X\longrightarrow X.$$

A measure $\mu$ defined in the space $X$ is called invariant with respect to $(\Phi(t,.))_{t\in \R}$ if for any $\mu-$measurable set $A$ one has
\begin{equation*}
\mu\big(\Phi(t,A)\big)=\mu(A), \quad t\in \R.
\end{equation*}
\end{defi}

In the case where $\mu$ is a probability measure,  the Poincar\'e theorem applies:

\begin{thm}[Poincar\'e theorem]
Let $(X, \mathcal{B}, \mu)$ be a probability space and let  $\Phi(t,.):X\longrightarrow X$ be a one parameter group which preserves  the  probability measure $\mu$.  
\begin{enumerate}[$(i)$]
\item Let $A\in \mathcal{B}$ be such that $\mu(A)>0$, then there exists $k\geq 1$ such that 
$$\mu\big(A \cap \Phi(k,A)\big)>0.$$
\item Let  $B\in \mathcal{B}$ be such that $\mu(B)>0$, then for $\mu$-almost all $x\in B$, the orbit $\big(\Phi(n,x)\big)_{n\in \N}$ enters infinitely many times in $B$.
\end{enumerate}
\end{thm}

 In the case of ordinary differential equations, the Liouville theorem provides a condition so that the Lebesgue measure (possibly with  a density) is invariant by the flow of the system. Namely, let $\Omega \subset \R^d$ be an open set and $F: \Omega \longrightarrow \R^d$ a $C^{\infty}$ function. Consider the ordinary differential equation 
 \begin{equation*}  
  \left\{
      \begin{aligned}
      & \dot{x}(t)=\frac{d x}{dt}(t)=F(x(t)),  \\[4pt]
       & x(0)=x_0.
      \end{aligned}
    \right.
\end{equation*}
We assume that for all $x_0\in \R$ the system has a unique solution $\Phi(t,x_0)$, such that ${\Phi(0,x_0)=x_0}$ and  which is defined for all $t\in \R$. The family $(\Phi(t,.)_{t\in \R}$ is a one parameter group of diffeomorphisms such that $\Phi(0,.)=id$, $\Phi(t,\Phi(s,.))=\Phi(t+s,.)$ for all $s,t\in \R$.

\begin{thm}[Liouville theorem]\label{thmLiou}
Denote by $dx$ the Lebesgue measure on $\Omega$ and let ${g :\Omega \longrightarrow [0,+\infty)}$ a $\mathcal{C}^{\infty}$ function. The flowmap $\Phi(t,.)$ preserves the measure $g dx$ if and only if
\begin{equation*} 
\sum_{k=1}^d\frac{\partial}{\partial x_k}\big(g F_k\big)=0.
\end{equation*}
\end{thm}
An important class of examples   is given by finite dimensional Hamiltonian systems of equations.


\subsection{Invariant measures for the Schr\"odinger equation on compact manifolds} 

  Let $M$ a compact manifold.  Then there exists a Hilbert basis $(h_n)_{n\geq 0}$ of $L^2(M)$, composed  of eigenfunctions of~$\Delta_M$ and we write 
$$-\Delta_M h_n= \lambda^2_n h_n\quad  \text{for all}\quad  n\geq 0.$$
 
Next, consider a  probability space
$(\Omega, {\mathcal F}, { \p})$ and let  $(g_n)_{n\geq 0}$  be a sequence of independent complex standard Gaussian variables $\mathcal{N}_{\C}(0,1)$
$$g_n=\frac1{\sqrt{2}}(g_{1,n}+i g_{2,n}),\qquad g_{1,n}, g_{2,n} \in \mathcal{N}_{\R}(0,1).$$

   Finally, let $(\alpha_n)_{n \geq 0}$ and define the probability measure $\mu$   via the map
 \begin{equation*} 
 \begin{array}{rcl}
\dis  \omega&\longmapsto &\dis\gamma({\om})=\sum_{n=0}^{+\infty} \alpha_n g_{n}(\om)h_n,  \qquad \mu= {\p} \circ \gamma^{-1}=\gamma_{\#} \p,
 \end{array}
 \end{equation*}
in other words: for all measurable set $A$, the measure $\mu$ is defined by
\begin{equation}\label{def-mu}
\mu(A)=\p\big(\om \in \Omega : \gamma(\om) \in A\big).
\end{equation}
 
 It is then easy to observe that   for any choice of  $(\alpha_n)_{n \geq 0}$ the measure $\mu$ is invariant by the flow of the linear Schr\"odinger equation:

 \begin{prop}\label{propinv}
 The measure $\mu$ defined in \eqref{def-mu} is invariant under the flow of the equation
  \begin{equation*}  
i\partial_sU+\Delta_MU=0,\quad (s,y) \in  \R\times M.
     \end{equation*}
\end{prop}

Considering the  nonlinear equation, it is then natural to look for invariant measures   (invariant Gibbs measures, see \cite{Bourgain2d} for example) or quasi-invariant measures (measures for which the nonlinear evolution is absolutely continuous with respect to $\mu$. In this direction, we refer to the works \cite{Tz2015,OhTz2,OhSoTz,OhTz1}).   \medskip

Let us sketch the proof of Proposition \ref{propinv}, since it is elementary.

\begin{proof} For all $t \in \R$, the random variable
$$\e^{it \Delta_M} \dis\gamma({\om})=\sum_{n=0}^{+\infty} \alpha_n e^{-it \lambda^2_n t}g_{n}(\om)h_n$$
 has the same distribution as $\gamma$ since for all $t \in \R$
 $$e^{-it \lambda^2_n t}g_{n}(\om) \sim \mathcal{N}_{\C}(0,1),$$
   hence the measure $\mu$ defined in \eqref{def-mu} is invariant by the linear flow.
 \end{proof}


  \subsection{The non compact case}
 Let us now turn to the case of the Schr\"odinger equation posed on~$\R$. Here the situation is dramatically different, since one has:
 
 \begin{prop}[\cite{Burq-Thomann}, Proposition 3.1]
Let $\s\in \R$ and consider a probability measure $\mu$ on $H^{\sigma}(\R)$. Assume that $\mu$ is invariant under the flow $\Sigma_{lin}$ of equation
 \begin{equation*}  
  \left\{
      \begin{aligned}
         &i\partial_sU+\partial_{y}^2U=0,\quad (s,y) \in  \R\times \R,
       \\  &  U(0,\cdot)  =U_0.
      \end{aligned}
    \right.
\end{equation*}
 Then $\mu=\delta_{0}$.
\end{prop}

By \cite[Proposition 3.2 and Proposition 3.3]{Burq-Thomann}, similar results  hold true for the nonlinear equation
$$i\partial_sU+\partial_{y}^2U=|U|^{p-1}U.$$

\begin{proof} Let us give the main lines of the argument. Let $\s\in \R$ and consider    $\mu$  an  invariant  probability measure on $H^{\sigma}(\R)$.  Let $\chi \in  {C}_0^{\infty}(\R)$. By invariance of the measure, 
\begin{equation*}  
 \int_{H^{\s}(\R)} \frac{\| \chi u \|_{H^{\sigma}}}{1+\|  u \|_{H^{\sigma}}}d \mu(u)= \int_{H^{\s}(\R)} \frac{\| \chi \Sigma_{lin}(t)u \|_{H^{\sigma}}}{1+\|  \Sigma_{lin}(t)u \|_{H^{\sigma}}}d \mu(u),
\end{equation*}
and by unitarity of the linear flow in $H^{\sigma}(\R)$, we get
\begin{equation}\label{chgt} 
 \int_{H^{\s}(\R)} \frac{\| \chi u \|_{H^{\sigma}}}{1+\|  u \|_{H^{\sigma}}}d \mu(u)= \int_{H^{\s}(\R)} \frac{\| \chi \Sigma_{lin}(t)u \|_{H^{\sigma}}}{1+\| u \|_{H^{\sigma}}}d \mu(u).
\end{equation}
\medskip
   
 Assume that the r.h.s.  of \eqref{chgt} tends to 0 when $t \to +\infty$. This implies  that $ \| \chi u \|_{H^{\sigma}}=0$, $\mu-$a.s., and thus  $\mu=\delta_0$ since $\chi$ is arbitrary.\medskip

  By continuity of the product by $\chi$ in  $H^{\sigma}(\R)$ and unitarity of the linear flow in~$H^{\sigma}(\R)$, we have
\begin{equation*} 
 \frac{ \| \chi \Sigma_{lin}(t)u \|_{H^{\sigma}}}{1+\| u \|_{H^{\sigma}}} \leq C  \frac{ \|   \Sigma_{lin}(t)u \|_{H^{\sigma}} }{1+\| u \|_{H^{\sigma}}}=C \frac{  \|    u \|_{H^{\sigma}}}{1+\| u \|_{H^{\sigma}}} \leq C.
\end{equation*} \medskip

     If $v \in \mathcal{C}_0^\infty(\R)$ is smooth,   by the  Leibniz rule and   dispersion
    \begin{equation*} 
 \| \chi \Sigma_{lin}(t)v  \|_{H^{\sigma}} \leq    \| \chi \|_{W^{\sigma,4}}     \|\Sigma_{lin}(t) v  \|_{W^{\sigma,4}}  \leq C t^{-1/4}  \|   v \|_{W^{\sigma,4/3}} \longrightarrow 0,
  \end{equation*}  
  when $t \longrightarrow +\infty$. We can conclude   with an approximation argument.
\end{proof}


  \subsection{Some functional analysis}  Define the harmonic oscillator
  $$H=-\partial^2_x+x^{2}\,.$$
There exists a Hilbert basis $(e_n)_{n\geq 0}$ of $L^2(\R)$, composed  of eigenfunctions of~$H$ and we write 
$$He_n=\lambda_n^2 e_n =  (2n+1) e_n\quad  \text{for all}\quad  n\geq 0.$$
  We define the harmonic Sobolev space $\W^{\s,p}$ by the norm $(\sigma>0)$
       \begin{equation*}  
      \Vert u\Vert_{\W^{\s,p}} = \Vert  H^{\s/2}u\Vert_{L^{p}} \equiv \Vert (-\Delta)^{\s/2} u\Vert_{L^{p}} + 
       \Vert\<x\>^{\s}u\Vert_{L^{p}}.
 \end{equation*}


  \subsection{Definition of the Gaussian measure $\mu_0$}

  Let $\eps>0$, we define the probability  Gaussian measure $\mu_0$ on $\H^{-\eps}(\R)$ as the distribution of the random variable $\gamma$
 \begin{equation*} 
 \begin{array}{rcl}
\Omega&\longrightarrow&\H^{-\eps}(\R)\\[3pt]
\dis  \omega&\longmapsto &\dis\gamma({\om})= \sum_{n=0}^{+\infty} \frac1{\lambda_n}g_{n}(\om)e_n, 
 \end{array}
 \qquad \mu_0= {\p} \circ \gamma^{-1}= \gamma _{\#} \p.
 \end{equation*}
 Notice that we can interpret $\mu_0$ as the  Gibbs measure  of the equation  $i \partial _t u-Hu = 0$. \medskip

   We   denote by 
$$X^0(\R)=  \bigcap_{\eps >0} \H^{-\eps}(\R).$$
 Thus  $\dis L^2(\R) \subset X^0(\R) \subset   \H^{-\eps}(\R)  $.
 One can show that  the measure $\mu_0$ satisfies $\mu_0(L^2(\R))=0$ and    $\dis \mu_0 \big( X^0(\R)\big)=1$, since $\mu_0(\H^{-\epsilon}(\R))=1$ for all $\epsilon>0$. This shows that the support of $\mu_0$ is essentially composed of $L^2$ functions. However, one can prove that the support of $\mu_0$  is actually smoother in  other $L^p$ scales. For instance, we have the bound:  
\begin{equation*} 
 \mu_0\big(\big\{u_0 \in X^0(\R) :   \|  \e^{-it H}u_0\|_{L^\infty{((-\pi, \pi)}; \W^{1/6-,\infty})}  \geq R\big\}\big)\leq C \e^{-cR^2},
 \end{equation*}
which implies that $\mu_0-$almost surely  $ \e^{-it H}u_0 \in L^\infty{((-\pi, \pi)}; \W^{1/6-,\infty})$.

   \subsection  {Equivalence of Gaussian measures}
  
Let $\mu$ and $\nu$ be two measures. We say that $\mu \ll \nu$ ($\mu$~is absolutely continuous with respect to $\nu$) if $\nu(A)=0 \implies  \mu(A)=0$. Now let us consider the particular case of Gaussian measures. Let   $\alpha_n, \beta_n>0$ and define  the measures $\mu={\bf p}\circ \gamma^{-1}$ and $\nu={\bf p}\circ \psi^{-1}$ with 
 \begin{equation*}
\dis\gamma({\om})= \sum_{n=0}^{+\infty} \frac1{\alpha_n}g_{n}(\om)e_n,\quad\quad \dis\psi({\om})= \sum_{n=0}^{+\infty} \frac1{\beta_n}g_{n}(\om)e_n.
\end{equation*}
 Then the measures $\mu$ and $\nu$ are absolutely continuous with respect to each other (this means that they have the same zero  measure sets) if and only if 
  \begin{equation}\label{condi}
  \sum_{n=0}^{+\infty}(\frac{\alpha_n}{\beta_n}-1)^2<+\infty.
  \end{equation}
 This criterion is equivalent to the next one which is written in a more symmetric manner by
   \begin{equation*}
  \sum_{n=0}^{+\infty}(\log \alpha_n -\log \beta_n)^2<+\infty,
  \end{equation*}
by an observation of R. Imekraz.  Criterion \eqref{condi} shows that it is very restrictive for two infinite-dimensional Gaussian measures to be absolutely continuous with respect to each other: the condition~\eqref{condi} says in some sense that the map which  $\gamma$ sends on $\psi$ has to be close to the identity. We refer to \cite[Appendix B.3]{Burq-Thomann} for more details.
  
 More generally,  in \cite[Section 2.2]{Burq-Thomann} we construct a four-parameter family of Gaussian measures based on the symmetries of the Schr\"odinger equation.


\section{Main results}
 
Consider the problem
  \begin{equation} \label{C1} 
  \left\{
      \begin{aligned}
         &i\partial_sU+\partial_{y}^2U=|U|^{p-1}U,\quad (s,y) \in  \R\times \R,
       \\  &  U(0)  =U_0 \in X^0(\R).
      \end{aligned}
    \right.
\end{equation}
\subsection{Global existence result}

We are now able to state our global existence result:

 \begin{thm}[\cite{Burq-Thomann}, Theorem 2.4]\label{thmglobal}
 Let $p>1$. 
  \begin{enumerate}[$(i)$]
\item    For $\mu_0-$almost every initial data $U_0 \in X^0(\R)$, there exists a unique, global in time, solution 
$$U=\Psi(s,0)U_0$$
 to~\eqref{C1}.  
\item   The  measures $\Psi(s,0)  _{{\#}} \mu_0 $  and $\Psi_{lin}(s, 0 )  _{{\#}} \mu_0$ are equivalent:
   \begin{equation} \label{star}
  \Psi_{lin}(s, 0 )  _{{\#}} \mu_0 \ll \Psi(s,0)  _{{\#}} \mu_0 \ll \Psi_{lin}(s, 0 )  _{{\#}} \mu_0.
\end{equation}
\item   For all $s' \neq s$, the measures $\Psi(s,0)  _{\#} \mu_{0}$ and  $\Psi(s',0)  _{\#} \mu_{0}$ are mutually singular.  
  \end{enumerate}
 \end{thm}
 It is interesting to compare this result with the case where NLS is posed on a compact manifold. In this latter case, the linear flow satisfies $\widetilde{\Psi}_{lin}(s, 0 )  _{{\#}} \mu_0=\mu_0$. Hence, in this context it is natural to try to prove a quasi-invariance result of the form $\widetilde{\Psi}(s,0)  _{{\#}} \mu_0 \ll  \mu_0$. But for NLS on the real line, we are able to prove, using criterion \eqref{condi} that $\Psi_{lin}(s, 0 )  _{{\#}} \mu_0$ and $\mu_0$ are mutually singular, hence $\Psi(s, 0 )  _{{\#}} \mu_0$ should not be compared to $\mu_0$ but to  $\Psi(s,0)  _{\#} \mu_{0}$ as in \eqref{star}.


\subsection{The scattering result}

 \begin{thm}[\cite{Burq-Thomann}, Theorem 2.4]\label{thm1} ~
 
  \begin{enumerate}[$(i)$]
  \item Assume that $p >1$. Then  for $\mu_0-$almost every initial data $U_0 \in X^0(\R)$, there exists  a  constant $C>0$ such that  for all  $s\in \R$
\begin{equation*}
\| \Psi(s,0)U_0\|_{L^{p+1}(\R)} \leq \begin{cases}
C\frac{(1+  \log\<s\> )^{1/(p+1)}} {\langle s\rangle ^{ \frac 1 2 -\frac 1 {p+1}}} &\text{ if } 1<p <5 \\[10pt]
\frac{C} {\langle s\rangle ^{ \frac 1 2 -\frac 1 {p+1}}} &\text{ if }  p\geq 5.
\end{cases}
\end{equation*}
\item Assume now that $p>3$. Then  there exist $  \eta>0$ and    ${W_\pm\in L^2(\R)}$ such that for all~$s\in \R$
   \begin{equation*}  
  \| \Psi(s,0) U_0- e^{is\partial_y^2} (U_0 +W_{\pm})\|_{L^2(\R)} \leq C \langle s \rangle^{- \eta}.
   \end{equation*}
     \end{enumerate}
\end{thm}

  For all $\phi \in \mathcal{C}_0^\infty(\R)$ we have the  dispersion bound
 \begin{equation*}
 \| \e^{is\partial_y^2}\phi\|_{L^{p+1}(\R)} \leq \frac{C} { |s|^{ \frac 1 2 -\frac 1 {p+1}}} \| \phi\|_{L^{(p+1)'}(\R)}, \qquad s \neq 0,
 \end{equation*}
therefore, the power decay in $s$ is optimal in the case $p \geq 5$ (we do not know if the $\log$ is necessary in the case $1<p<5$). We also stress that the condition $p>3$ is optimal in our scattering result, by \cite{Barab}.\medskip

  We conclude this section by giving a few references on scattering results for NLS. \medskip

$\bullet$ Deterministic scattering results for NLS: \medskip

 \begin{tabular}{ll}
Barab \cite{Barab}&$\longrightarrow$ never scattering when $p\leq 3$ \\
Tsutumi--Yajima \cite{TsYa} \quad\quad&  $\longrightarrow$ scattering in $L^2$ with $H^1$ data  \\
Nakanishi \cite{Kenji} &   $\longrightarrow$ scattering in $H^\sigma$    \\
Dodson \cite{Dodson} &  $\longrightarrow$   scattering in $L^2$    when $p = 5$     \\
\end{tabular}
 \medskip
 
$\bullet$ Probabilistic scattering results for NLS: \medskip

 \begin{tabular}{ll}
Burq--Thomann--Tzvetkov \cite{BTT}&$\longrightarrow$ case $d=1$ and $p\geq 5$ \\
Poiret--Robert--Thomann \cite{PRT2} & $\longrightarrow$ case $d\geq 2$ and $p\geq 3$  \\
Dodson--L\"uhrmann--Mendelson \cite{DLM} \quad \quad & $\longrightarrow$ case $d=4$ and $p=3$  \\
Killip--Murphy--Visan \cite{KMV} &   $\longrightarrow$ case $d=4$  in the radial setting    \\
Latocca \cite{Latocca} &  $\longrightarrow$ case $d\geq 2$ in the radial setting      \\
Camps \cite{Camps} & $\longrightarrow$ case $d=3$ and $p= 3$  in the radial setting   \\
\end{tabular}
\medskip

 For results on the Gross-Pitaevskii equation with random perturbations we refer to the works of de Bouard, Debussche and Fukuizumi \cite{BDF1, BDF2, BDF3} and references therein.

 \section{Some key ingredients of the proof}
 
\subsection{The lens transform: compactification in time and space}
 
There is an explicit transform, called the lens transform, which maps the solutions of NLS to solutions of NLS with harmonic potential. Namely,  if $U(s,y)$ is a solution of the problem \eqref{C1}, then the function  $u(t,x)$ defined  for $|t|<\frac{\pi}{4}$ and $x\in\R$ by
\begin{equation*} 
 u(t,x)=  \mathscr{L} (U)(t,x):= \frac{1}{\cos^{\frac{1}{2}}(2t)} U\big(\frac{\tan(2t)} 2, \frac{x}{\cos(2t)}\big)e^{-\frac{ix^2{\rm tan}(2t)}{2}}
 \end{equation*}
 solves the problem
 \begin{equation}\label{NLSp}
  \left\{
      \begin{aligned}
         &i\partial_t u-Hu=\cos^{\frac{p-5}{2}}(2t)|u|^{p-1}u,\quad |t|<\frac{\pi}{4},\, x\in\R,
       \\  &  u(0,\cdot)  =U_0.
      \end{aligned}
    \right.
\end{equation}

We define the corresponding energy 
\begin{equation}\label{energ}
\mathcal{E}(t,u(t))=\frac{1}{2}\|\sqrt{H}\,u(t)\|_{L^2(\R)}^2+\frac{\cos^{\frac{p-5}{2}}(2t)}{p+1}
\|u(t)\|_{L^{p+1}(\R)}^{p+1}\, ,
\end{equation}
which  is  not conserved.   For $-\frac{\pi}4<t<\frac{\pi}4$, we define the measure 
 \begin{equation*}
 d\nu_{t}= e^{-\mathcal{E}(t,u)}dud\ov{u}  = e^{-\frac{\cos^{\frac{p-5}{2}}(2t)}{p+1}\|u\|_{L^{p+1}(\R)}^{p+1}}d\mu_0 
  \end{equation*}
which is therefore not invariant by the flow of \eqref{NLSp}.


\subsection{Monotonicity of the measure $\nu_t$} 
 We are able to bound the nonlinear evolution of $\nu_0$ by~$\nu_t$. More precisely, we have:
  \begin{prop}\label{prop-niko}
  Denote by $\Phi(t,0)$ the flow of equation \eqref{NLSp}. Then, for all $0\leq |t| <\frac{\pi}4$
  
   \begin{align} \label{mono}
\nu_{0}(\Phi(t,0) ^{-1}A )&\leq \begin{cases}\big[\nu_{t} (A)\big] ^{\smash{\big( {\cos(2t)} \big) ^{\frac{ 5-p} 2}}}&\text{ if } 1\leq p \leq 5\\[5pt]
 \nu_{t} (  A) &\text{ if }  p \geq 5. \end{cases}
 \end{align}
   \end{prop}
  
These quantitative estimates will be in the core of our argument. In particular, they allow to extend the globalisation argument of Bourgain relying on invariant measures  (see Section~\ref{sect-Bourgain}). They are also crucial in the proof of the scattering result. The monotonicity estimate \eqref{mono} has already been used in \cite{BTT} in the case $p \geq 5$ in order to prove global well posedness results for NLS. Notice that in the case $1 \leq p \leq 5$, the estimates \eqref{mono} are quite accurate for small times, but they deteriorate when $|t|$ is close to $\pi/4$: in this regime we recover the trivial bound $\nu_{0}(\Phi(t,0) ^{-1}A ) \leq 1$.
\medskip

We give here an outline of the proof of Proposition \ref{prop-niko}.  Recall the definition \eqref{energ}, then a direct computation shows that 
\begin{equation*} 
\frac{d}{dt}\big(\mathcal{E}(t,u(t))\big)=
\frac{(5-p)\sin(2t)\cos^{\frac{p-7}{2}}(2t)}{p+1}
\|u(t)\|_{L^{p+1}(\R)}^{p+1}\,.
\end{equation*}
Next, set $F(t)=\nu_{t}\big(\Phi(t,0)A\big)$. Then we have to prove that  for all $0\leq |t| <\frac{\pi}4$
   \begin{align*} 
F(0)&\leq \begin{cases}\big[F(t)\big] ^{\smash{\big( {\cos(2t)} \big) ^{\frac{ 5-p} 2}}}&\text{ if } 1\leq p \leq 5\\[5pt]
F(t)&\text{ if }  p \geq 5. \end{cases}
 \end{align*}
 We compute 
 $$
 \frac{d}{dt} F(t) =   (p-5) \tan(2t)  \int_A \alpha\big(t,u(t)\big)  e^{-\mathcal{E}(t,u(t))}du_0,
 $$
where $\alpha(t,u)=\frac{\cos^{\frac{p-5}{2}}(2t)}{p+1}\|u\|_{L^{p+1}(\R)}^{p+1}$. Using the H\"older inequality, we can check that for all $k \geq 1$
 \begin{equation*}
\frac{d}{dt}F(t) \leq  (p-5) \tan(2t) \frac { k} e\big( F(t)\big)^{1-\frac1k}.
\end{equation*}\medskip
Next, optimizing  with $k = - \log\big( F(t)\big)$ yields
   $$ 
\frac{d}{dt}F(t) \leq    - (p-5) \tan(2t) \log\big( F(t)\big) F(t).
$$ 
Finally the result follows from the integration of the previous differential inequality.


 \subsection{On Radon-Nikodym derivatives}
The bounds obtained in Proposition~\ref{prop-niko} say much more than just an absolute continuity result between two measures. In fact, they provide integrability results on the    Radon-Nikodym density, since one has the following general result:

 \begin{prop}[\cite{Burq-Thomann}, Proposition 3.5]\label{prop-R-N}
Let $\mu, \nu$ be two finite  measures on a measurable space~$(X,\mathcal{T})$. Assume that 
\begin{equation*} 
\mu \ll \nu,
\end{equation*} and more precisely 
\begin{equation} \label{R-N}
\exists \,  0< \alpha \leq1, \quad \exists \,C>0, \quad \forall A \in \mathcal{T}, \quad \mu(A) \leq C \nu(A)^{\alpha}.
\end{equation}
By the  Radon-Nikodym theorem,  there exists a $f\in L^1(d\nu)$  with $f\geq 0$, such that  $ \dis d \mu = f d\nu$, and we write   $\dis f= \frac{d\mu}{d\nu}$.
 \begin{enumerate}[$(i)$]
\item The assertion~\eqref{R-N} is satisfied with $0< \alpha<1$ iff  $f\in L^p_{w}(d\nu)\cap L^1(d\nu)$ with   $p = \frac 1 {1- \alpha}$. In other words, $f\in L^1(d\nu)$ and 
$$\nu \big(\big\{ x : \; |f(x) |\geq \lambda \big\}\big) \leq C' \langle \lambda\rangle ^{- p}, \qquad \forall\, \lambda >0.$$
\item The assertion~\eqref{R-N} is satisfied with $ \alpha=1$ iff  $f\in L^{\infty}(d\nu)\cap L^1(d\nu)$.
\end{enumerate}
\end{prop}


\subsection{The Bourgain argument revisited}\label{sect-Bourgain}
 
 Let us now show how local in time solutions can be extended, using the bound \eqref{mono}. In order to simplify the presentation of the following argument, we assume that $p \geq 5$. Moreover, we do not give details on the norm $\|\;\cdot \;\|$ below, since it does not really play a role in the method. Thus, let us assume the three following facts: \medskip
 
  $\blacktriangleright$ There exists a flow $ \Phi $ such that the time of existence $ \tau $ on the ball
\begin {equation*}
B_{R} = \big \{u \in X ^ {0} (\R) \;: \; \| u \| \leq R ^ {1/2} \, \big \},
\end {equation*}
is uniform and such that $ \tau \sim R ^ {- \kappa} $ for some $ \kappa> 0 $.  \medskip\medskip

  $\blacktriangleright$ For all $ | t | \leq \tau $
\begin{equation*} 
\Phi (t,0) \big (B_{R} \big) \subset \big \{u \in X ^ {0} (\R) \;: \; \| u \| \leq (R + 1) ^ {1/2} \, \big \}.
\end {equation*}
  \medskip

  $\blacktriangleright$  We have the large deviation estimate $ \mu_0 (X ^ {0} (\R) \backslash B_{R}) \leq C \e ^ {- cR} $.
 \medskip

Then for $ T \leq \e ^ {cR / 2} $ fixed, we define the set of the good data
\begin{equation*} 
\Sigma_{R} = \bigcap_{k = - [T / \tau]} ^ {[T / \tau]} \Phi (k \tau,0)^{-1} \big (B_{R} \big).
\end{equation*} 
By Proposition \ref{prop-niko} we have
\begin{eqnarray*}
\nu_0 (X ^ {0} (\R) \backslash \Sigma_{R}) & \leq & \sum_{k = - [T / \tau]} ^ {[T / \tau]}   \nu_0 \Big(  \Phi (k \tau,0)^{-1}\big(   X ^ {0} (\R) \backslash B_{ R}\big)\Big)  \\
& \leq & \sum_{k = - [T / \tau]} ^ {[T / \tau]}   \nu_{k \tau}    \big(   X ^ {0} (\R) \backslash B_{ R}\big).
\end{eqnarray*}
In his  original  argument, Bourgain \cite{Bourgain2d} considered invariant measures, so that the previous estimate was an indeed equality in his case.   We observe here that the monotonicity property \eqref{mono} is sufficient. Next, by definition of $\nu_t$, we have $\nu_t(A) \leq \mu_0(A)$, so that   
\begin{eqnarray*}
\nu_0 (X ^ {0} (\R) \backslash \Sigma_{R}) &\leq &   \big(2 [T / \tau]+1\big)   \mu_0\big(   X ^ {0} (\R) \backslash B_{ R}\big)    \\
& \leq &  c\e ^ {-cR / 2}\end{eqnarray*}
which shows that $ \Sigma_{R} $ is a big set of $ X ^ {0} (\R) $ when $ R \longrightarrow + \infty $. \medskip

We deduce that for all $ | t | \leq T $ and $ u \in \Sigma_{R} $
\begin {equation*}
\| \Phi (t,0) (u) \| \leq (R + 1) ^ {1/2}.
\end{equation*}
In particular, for $ | t | = T \sim \e ^ {cR / 2} $
\begin{equation*}
\| \Phi (t,0) (u) \| \leq C (\ln | t | +1) ^ {1/2}.
\end{equation*}


\end{document}